\documentclass{amsart}
\usepackage{amsfonts,amssymb,epsfig}
\usepackage[latin1]{inputenc}
\usepackage[all]{xy}
\usepackage{amsmath}
\usepackage{delarray}

\newtheorem{thm}{Theorem}[section]
\newtheorem{cor}[thm]{Corollary}
\newtheorem{lem}[thm]{Lemma}
\newtheorem{prop}[thm]{Proposition}
\theoremstyle{definition}
\newtheorem{defn}[thm]{Definition}
\theoremstyle{remark}
\newtheorem{rem}[thm]{Remark}
\numberwithin{equation}{section}
\newtheorem{example}[thm]{Example}

\newcommand{\ove}{\overline}

\newcommand{\N}{{\mathbb N}}

\newcommand{\R}{{\mathbb R}}
\newcommand{\C}{{\mathbb C}}
\newcommand{\Proj}{{\mathbb P}}
\newcommand{\X}{\mathfrak{X}}

\def\lcf{\lbrack\! \lbrack}
\def\rcf{\rbrack\! \rbrack}

\begin{document}

\title[A cohomology attached to a function]
{A cohomology attached to a function}

\author{Philippe Monnier}
\address{Departamento de Matem\'atica, Instituto Superior T\'ecnico,
Avenida Rovisco Pais, 1049-001 Lisbon, Portugal} 
\email{pmonnier@math.ist.utl.pt}
\keywords{Differential forms, Poisson cohomology, Lie algebroids,
  singularities}
\subjclass{53D17, 58A10, 14F40, 32S05 }

\begin{abstract}
In this paper, we study a certain cohomology attached to a smooth
function, which arose naturally in Poisson geometry. We explain how
this cohomology depends on the function, and we prove that it
satisfies both the excision and the Mayer-Vietoris axioms. For a
regular function we show that the cohomology is related to the de Rham
cohomology. Finally, we use it to give a new proof of a well-known
result of A. Dimca in complex analytic geometry.
\end{abstract}

\maketitle

\section{Introduction}
There are several cohomologies attached to a function that can be
defined in terms of differential forms, such as the relative cohomology
associated to a singularity, or the cohomology of the complex of 
logarithmic differential forms associated with the complement of an
hyperplane. These cohomologies give, for instance, information on the
topology of the complement of the zeros of the function. In this
paper, we consider a new cohomology attached to a smooth function on a
differentiable manifold. 

This new cohomology is also defined in terms of differential
forms. More precisely, if $M$ is a differentiable manifold and $f$ is
a smooth function on $M$, we define a coboundary operator
\begin{align*}
d_f:\Omega^k(M) &\longrightarrow \Omega^{k+1}(M)\\
         \alpha &\longmapsto  fd\alpha-kdf\wedge\alpha.
\end{align*}
where $\Omega^k(M)$ is the space of $k$-differential forms on $M$.
It is easy to check that $d_f\circ d_f=0$, and we denote by $H_f^\bullet(M)$ 
the cohomology associated with the complex $(\Omega^\bullet(M),d_f)$.
More generally, for any integer $p$, we define a coboundary operator
\begin{align*}
d_f^{(p)}:\Omega^k(M) &\longrightarrow \Omega^{k+1}(M)\\
               \alpha &\longmapsto  fd\alpha-(k-p)df\wedge\alpha.
\end{align*}
We still have $d_f^{(p)}\circ d_f^{(p)}=0$ and we denote by
$H_{f,p}^\bullet(M)$ the cohomology of this complex. We shall restrict
our attention to the cohomology $H_f^\bullet(M)$ but most results
readily generalize to the cohomology $H_{f,p}^\bullet(M)$.

This cohomology was considered for the first time in \cite{M2} in
the context of Poisson geometry, and more generally, Nambu-Poisson
geometry. There we have computed this cohomology in the case where $f$
is the germ of a function with an isolated singularity. The aim of
this paper is to initiate a systematic study of this cohomology. 

We start, in Section \ref{section:origins}, by showing several
possible ways of defining this cohomology. First we recall how it
arises in Poisson and Nambu-Poisson geometry.  Then we construct a
certain Lie algebroid attached to a function $f$ for which the Lie
algebroid cohomology coincides with $H_f^\bullet(M)$. For a regular
function $f$, there is another Lie algebroid one can attach to the
function, namely the Melrose fake tangent bundle of $S=f^{-1}(\{0\})$
(see \cite{CW}). This Lie algebroid does not coincide with ours, but
they have isomorphic cohomologies. Finally, one can also consider
differential forms with a ``pole'' along $S$, obtaining a chain complex for
which the cohomology is also $H_f^\bullet(M)$.

In Section \ref{section:properties} we study some basic properties of the
cohomology. First we discuss how the cohomology varies when the
function $f$ changes. In particular, we show that if the function $f$
does not vanish, then the cohomology $H_f^\bullet(M)$ coincides
with the de Rham cohomology of $M$. Then we will show that it is
possible to write a Mayer-Vietoris exact sequence, a relative
cohomology exact sequence, and an excision theorem, for our
cohomology. We also give an appropriate notion of homotopy, but it is
an open question whether the cohomology is homotopy invariant in
general.

%The problem arised from the homotopy invariance question.
%We give some answers in the regular case (in section 6) but in the general 
%case, we do not know yet if this cohomology is homotopy invariant.

In Section \ref{section:regular} we consider the regular case, i.e.,
the case where the function $f$ does not have singularities in a
neighborhood of $S=f^{-1}(\{0\})$. In this case, we can relate the
cohomology with the de Rham cohomology of $M$ and of $S$, showing that
the space $H^k_f(M)$ is isomorphic to $H^k_{dR}(M)\oplus
H^{k-1}_{dR}(S)$. As a corollary of this result, one obtains
the Poisson cohomology for generic 2-dimensional Poisson
structures. In the regular case, we prove homotopy invariance.

Finally, in the last section, we study the complex case, giving an
application of our cohomology to complex algebraic geometry. Namely, 
we explain how the results we have found in \cite{M2} can be applied
to give information on the degeneration of a spectral sequence
converging to the cohomology of an hypersurface complement. As a
corollary, we obtain a new proof of a well-known result of A. Dimca.

\section*{Acknowledgments}
I would like to thank Jean-Paul Dufour for his help and its comments. 
I wish also to thank Alexandru Dimca, Rui Loja Fernandes and Alan Weinstein for their remarks.
%%%%%%%%%%%%%%%%%%%%%%%%%%%%%%%%%%%%%%%%%%%%%%%%%%%%%%%%%%%%%%%%%%%%%%%%%%%%%%%

\section{Geometrical origins}
\label{section:origins}
The following notations will be enforced throughout the paper. We
denote by $f$ a smooth function on a $n$-dimensional manifold $M$ and
by $S\subset M$ the level set $f^{-1}(\{0\})$. As usual, $\Omega^k(M)$
denotes the vector space of $k$-differential forms, and $H_{dR}^k(M)$
the $k$-th de Rham cohomology group.  Dually, $\X^k(M)$ denotes the
vector space of $k$-vector fields. Also,
$[~,~]:\X^k(M)\times\X^l(M)\to\X^{k+l-1}(M)$ denotes the Schouten
bracket on multi-vector fields. For a cohomology theory, we denote by
$Z^k$ (resp. $B^k$) the space of $k$-cocycles (resp. $k$-cobords).

%%%%%%%%%%%%%%%%%%%%%%%%%%%%
\subsection{The two-dimensional case}
\label{sec:2 dim}
Let $M$ be a Poisson manifold with Poisson 2-vector field
$\Pi\in\X^2(M)$, so that $[\Pi,\Pi]=0$ (see for instance
\cite{CW,L,V}). If the manifold $M$ has dimension two, this condition
is automatically satisfied, so every 2-vector on a 2-dimensional
manifold is a Poisson structure.

Assume that $(M,\Pi)$ is a 2-dimensional orientable Poisson manifold,
and fix a volume form $\nu\in\Omega^2(M)$. The contraction $f:=i_\Pi\nu$
is a smooth function. We have observed in \cite{M2} that the \emph{Poisson
cohomology} of $(M,\Pi)$ is isomorphic to $H_f^\bullet(M)$. Let us
recall how this works.

First of all, the Poisson cohomology of $(M,\Pi)$ is defined to be the
cohomology of the following chain complex (see \cite{L}):
\[
\xymatrix{
0\ar[r]&\X^0(M)\ar[r]^{\partial}&\X^1(M)\ar[r]^{\partial}&\X^2(M)\ar[r]&0
}
\]
where the boundary map is $\partial(Q)=[Q,\Pi]$. Hence, the map
$\partial:\X^0(M)\to \X^1(M)$ is the map that associates to a function
$g$ its Hamiltonian vector field $X_g$:
\[ \partial (g)=[g,\Pi]\equiv X_g,\]
and $\partial:\X^1(M)\to \X^2(M)$ is the map that associates to a
vector field $X$ the Lie derivative of $\Pi$ along $X$:
\[ \partial (X)=[X,\Pi]\equiv \mathcal{L}_X\Pi.\]
This cohomology is an invariant of the Poisson manifold, which 
has been studied, from different points of view, for instance in
\cite{M1,N,R,RV,V}.

Secondly, we have an isomorphism of chain complexes 
\[
\phi:(\X^\bullet(M),\partial)\longrightarrow (\Omega^\bullet(M),d_f),
\]
where $\phi^0:C^\infty(M)\to C^\infty(M)$ is the identity,
$\phi^1:\X^1(M)\to \Omega^1(M)$ is contraction of $\nu$:
\[\phi^1(X)\equiv-i_X\nu,\]
and $\phi^2:\X^2(M)\to \Omega^2(M)$ is the linear application
defined by
\[\phi^2(\Gamma)\equiv(i_\Gamma\nu)\nu.\]

The Poisson cohomology of a manifold is, in general, very hard to
compute, even in dimension two. Since working with differential
forms has many advantages over working with multivectors, one may expect that
this isomorphism will lead to actual computations of Poisson cohomology
in dimension two. We shall see an example of that in the proof of 
Theorem \ref{thm:Poisson cohomology 2}.

%%%%%%%%%%%%%%%%%%%%%%%%%%%%%
\subsection{In higher dimensions}
If $M$ is an orientable manifold of dimension $n>2$, one generalizes
the previous case in a straightforward way. One considers a $n$-vector
$\Lambda\in\X^n(M)$, and fixes a volume form $\nu\in\Omega^n(M)$,
obtaining a smooth function $f:=i_\Lambda\nu$. The pair $(M,\Lambda)$
is no more a Poisson manifold, but it is a \emph{Nambu-Poisson
manifold} of degree $n$, which may be seen as a kind of generalization
of Poisson structures (see \cite{Na,T}).

Now we would like to associate a cohomology to the pair $(M,\Lambda)$,
generalizing Poisson cohomology in dimension two. In \cite{I}, the
authors construct a chain complex (called the Nambu-Poisson complex)
associated to any Nambu-Poisson manifold of dimension and of degree
larger than 3. This complex is rather difficult to manipulate, but we
have shown in \cite{M2} that the Nambu-Poisson cohomology of
$(M,\Pi)$ is indeed isomorphic to $H_f^\bullet(M)$.

There is a second complex one can associate to the pair $(M,\Lambda)$,
which also generalizes Poisson cohomology in dimension two, and which
is much simpler. One takes
\[
\xymatrix{
0\ar[r]&\big(C^\infty(M)\big)^{n-1} \ar[r]^{\partial}&
\X^1(M)\ar[r]^{\partial}&\X^n(M)\ar[r]&0
}
\]
where the boundary map $\partial:\X^0(M)\to \X^1(M)$ is the map that
associates to the functions $g_1,\dots,g_n$ their Hamiltonian vector
field $X_{g_1,\dots,g_{n-1}}$: 
\[ \partial (g_1,\dots,g_{n-1})=i_{dg_1\wedge\hdots\wedge dg_{n-1}}
\Lambda\equiv X_{g_1,\dots,g_{n-1}},\]  
and $\partial:\X^1(M)\to \X^n(M)$ is the map that associates to a
vector field $X$ the Lie derivative of $\Lambda$ along $X$:
\[ \partial (X)=[X,\Lambda]\equiv \mathcal{L}_X\Lambda.\]
In the same way as for the 2-dimensional case, one can show that the
last two cohomology groups of this chain complex are isomorphic to 
$H^{n-1}_{f,n-2}(M)$ and $H^n_{f,n-2}(M)$ (see \cite{M2}).

%%%%%%%%%%%%%%%%%%%%%%%%%%%%%%
\subsection{A Lie algebroid attached to a function}
\label{sec:Lie algebroid}
Recall (see, e.g., \cite{CW,F,Mk}) that a \emph{Lie algebroid} over
$M$ is a triple $(A,\rho,\lcf\,,\,\rcf)$ where $A$ is a vector bundle
over $M$, $\rho:A\to TM$ is a bundle map (called the \emph{anchor}),
and $\lcf\,,\rcf$ is a Lie algebra bracket on the sections
$\Gamma(A)$, such that:
\begin{itemize}
\item $\rho$ defines a Lie algebra homomorphism
  $(\Gamma(A),\lcf\,,\,\rcf)\to(\X(M),[~,~])$; 
\item for every $u,v\in\Gamma(A)$ and $g\in C^\infty(M)$:
\[ \lcf u,gv\rcf=g\lcf u,v\rcf+\big( \rho(u)\cdot g\big)v.\]
\end{itemize}
To any Lie algebroid one associates a cohomology $H^\bullet(A)$ 
by considering the chain complex $(\Omega^\bullet(A),d_A)$, where 
$\Omega^k(A)\equiv\Gamma(\wedge^k A^*)$ and 
\begin{multline*}
d_A Q(u_0,\dots,u_r)=\frac{1}{r+1}\sum_{k=0}^{r}
(-1)^k\rho(u_k)\cdot Q(u_0,\dots,\widehat{u}_k,\dots,u_r)\\
+\frac{1}{r+1}\sum_{k<l}(-1)^{k+l+1}Q([u_k,u_l],u_0,\dots,\widehat{u}_k,\dots,\widehat{u}_l,\dots,u_r).
\end{multline*}

Now, for any smooth function $f$ on a manifold $M$ we can attach a Lie
algebroid as follows. We take $A=TM$, the anchor $\rho:TM\to TM$ is
defined by 
\[ \rho(X)\equiv fX,\quad X\in\X(M),\] 
and the Lie bracket $\lcf\,,\,\rcf$ on $\X(M)$ is given by
\[ 
\lcf X,Y\rcf\equiv\frac{[fX,fY]}{f}=f[X,Y]+(X\cdot f)Y-(Y\cdot f)X,
\quad X,Y\in\X(M).
\]
It is easy to check that the triple $(TM,\rho,\lcf\,,\,\rcf)$ is a Lie
algebroid over $M$ and its cohomology is precisely $H^\bullet_f(M)$.

\begin{rem}
  The Lie algebroid $(TM,\rho,\lcf\,,\,\rcf)$ is always integrable to
  a Lie groupoid since the obstructions to integrability given in
  \cite{CF} vanish. 
%  Adopting the same notations as in \cite{CF}, for
%  each $x$ in $M$ we let $L_x$ be the leaf of the algebroid through
%  $x$ and we let $\mathfrak{g}_x\equiv\ker\rho_x$ be the isotropy Lie
%  algebra. We also denote by $G(\mathfrak{g}_x)$ the simply-connected
%  Lie group integrating $\mathfrak{g}_x$. In \cite{CF}, it is defined
%  a certain group homomorphism $\partial:\pi_2(L_x,x)\to
%  G(\mathfrak{g}_x)$ with image $N_x$ (called the \emph{monodromy
%  group} at $x$) giving the obstructions to the integrability of the
%  algebroid (\cite{CF}, Theorem 4.1). For the Lie algebroid attached
%  to a function, we have two possibilities: (i) points $x\in M$ where
%  $f(x)=0$, in which case the leaf $L_x$ reduces to $\{x\}$, so
%  $\pi_2(L_x,x)$ is trivial; and (ii) points $x\in M$ where $f(x)\neq
%  0$, in which case $\rho_x$ is injective, so$\mathfrak{g}_x$ is
%  trivial. In any case, $N_x$ is trivial so the obstructions vanish.
\end{rem}

\begin{rem} 
  When the function $f$ is regular there is another Lie
  algebroid attached to $f$ which can be defined as follows (see
  \cite{Me} and \cite{CW}). Recall that $S\subset M$ denotes the set
  $f^{-1}(0)$, which is an embedded submanifold if $f$ is regular. It
  is shown in \cite{Me}, that the $C^\infty(M)$-module $\X_S(M)$ of
  vector fields on $M$ tangent to $S$ is the space of sections of a
  vector bundle $A$ over $M$, called the \emph{fake tangent
  bundle}. On $A$ one has a structure of a Lie algebroid over $M$,
  where the bracket is the standard Lie bracket of vector fields, and
  the anchor may be defined locally as follows. For a point $p\in S$,
  there exists local coordinates $(U,x,y_2,\hdots,y_n)$ such that
  $U\cap S=\{q\in U:x(q)=0\}$. If one sets $e_1=x\frac{\partial
  }{\partial x}$ and $e_i=\frac{\partial}{\partial y_i}$ for $i>1$,
  the $e_i$'s form a local basis of $\X_S(M)$. The anchor map $\tau$
  is then defined as $\tau(e_1)=x\frac{\partial}{\partial x}$ and
  $\tau(e_i)=\frac{\partial}{\partial y_i}$ for $i>1$. This Lie
  algebroid does not coincide with the one defined above (the later
  has points of rank zero, while the first one not), but we will see
  later (cf.~Remark \ref{rk:cohomelrose}) that their Lie algebroid
  cohomologies are isomorphic.
\end{rem}

\begin{rem} 
  {For} $p\neq 0$ the operator $d_f^{(p)}$ is not a
  derivation of the exterior algebra, hence the cohomology
  $H_{f,p}^\bullet(M)$ does not come from a Lie algebroid. 
\end{rem}

%%%%%%%%%%%%%%%%%%%%%%%%%%%%%%%%%%%%%%%%%%%%%%%%%%%%%%%%%%%%%%%%%%%%%%%%%%%%%%%
\subsection{Singular k-forms} 
Let us call a form $\omega\in\Omega^k(M\setminus S)$ a \emph{singular
$k$-form} if the form $f^k\omega$ can be extended to a smooth form on
$M$. We denote the space of singular $k$-forms by $\Omega_f^k(M)$.

If $\omega\in\Omega_f^k(M)$ is a singular $k$-form then $d\omega$ is a
singular $(k+1)$-form. In fact, we have
\[f^{k+1}d\omega=d(f^{k+1}\omega)-(k+1)df\wedge (f^k\omega),\]
so $f^{k+1}d\omega$ also extends to a smooth form on $M$.  Therefore we
obtain a chain complex $\big(\Omega_f^\bullet(M),d\big)$.

\begin{prop}
\label{autreexpress}
The cohomology of $\big(\Omega_f^\bullet(M),d\big)$ is isomorphic to
$H_f^\bullet(M)$. 
\end{prop}

\begin{proof}
Define a map of chain complexes
$\varphi:\big(\Omega_f^\bullet(M),d\big)\to
\big(\Omega^\bullet(M),d_f\big)$
by setting
\[
\varphi^k:\Omega_f^k(M)\to \Omega^k(M),\quad \omega \mapsto f^k\omega.
\]
It is easy to check that $\varphi$ induces an isomorphism in cohomology.
%If $\omega\in\Omega^k_f(M)$ satifies $\varphi^k(\omega)=0$, then 
%$f^k\omega=0$ on $M$ and so $\omega=0$ on $M\setminus S$. We deduce that
%$\phi^k$ est injective.\\
%Now, let $\omega$ be in $\Omega^k(M)$. We consider the $k$-form $\alpha$ 
%on $M\setminus S$ defined by $\alpha=\frac{1}{f^k}\omega$. 
%We then have $\alpha\in\Omega_f^k(M)$ and 
%$\phi^k(\alpha)=\omega$. Therefore, $\phi^k$ is surjective.
\end{proof}

%%%%%%%%%%%%%%%%%%%%%%%%%%%%%%%%%%%%%%%%%%%%%%%%%%%%%%%%%%%%%%%%%%%%%%%%%%%%%%%
\section{Basic Properties}
\label{section:properties}
In this section we will study some basic properties of the
cohomology defined above.

\subsection{Degree zero cohomology} If $M\setminus S$ is a dense
subset of $M$ (e.g., if $f$ is regular) one can compute the groups
$H^0_{f,p}(M)$:

\begin{prop}
\label{prop:H0cst}
If $M\setminus S$ is dense in $M$,
\[
H^0_{f,p}(M)=
\left\{
\begin{array}{l}
0,\quad \text{ if }p>0,\\
\\
\R,\quad \text{ if }p\le 0.
\end{array}
\right.
\]
\end{prop}

\begin{proof}
If $p>0$, note that $d_f^{(p)}(g)=\frac{d(f^pg)}{f^{p-1}}$ for any
smooth function $g$ on $M$. Hence $d_f^{(p)}(g)=0$ iff $g\equiv 0$,
and we obtain $H^0_{f,p}(M)=\{0\}$.

If $p\le 0$, let $g$ be a function on $M$ such that $d_f^{(p)}(g)=0$. We have 
$d\big( \frac{g}{f^{-p}}\big)=0$ on $M\setminus S$, so $g=\lambda
f^{-p}$ on $M\setminus S$ for some $\lambda\in\R$. It follows that $g=\lambda
f^{-p}$ on $M$, so we obtain $H^0_{f,p}(M)\simeq\R$.
\end{proof}

The higher degree cohomology groups are much harder to compute, even
in the case where the function vanishes at a single point.

\subsection{Dependence on the function}
\label{section:dependence}
A natural question to ask about the cohomology $H_f^\bullet(M)$ is
how it depends on the function $f$.  A first result is the following.

\begin{prop}
If $h\in C^\infty(M)$ does not vanish, then the cohomologies
$H_f^\bullet(M)$ and $H_{fh}^\bullet(M)$ are isomorphic.
\end{prop}

\begin{proof}
For each $k\in\N$, consider the linear isomorphism 
\[
\phi^k :\Omega^k(M)\to\Omega^k(M),\qquad 
\alpha\longmapsto\frac{\alpha}{h^k}.
\]
If $\alpha$ is a $k$-form on $M$, one checks easily that
\[ 
\phi^{k+1}(d_{fh}\alpha)=d_f\big(\phi^k(\alpha)\big),
\]
so $\phi$ induces an isomorphism between the cohomologies
$H_f^\bullet(M)$ and $H_{fh}^\bullet(M)$.
\end{proof}

\begin{cor}
If the function $f\in C^\infty(M)$ does not vanish, then
$H_f^\bullet(M)$ is isomorphic to the de Rham cohomology
$H_{dR}^\bullet(M)$.
\end{cor}

It follows also that the cohomology $H_f^\bullet(M)$ depends only on
the germ of the function $f$ on its set of zeros:

\begin{cor}
\label{cor:deformation}
If $g$ and $f$ are smooth functions on $M$ such that
$S=f^{-1}(0)=g^{-1}(0)$ and $g=f$ on some neighborhood of $S$, then
$H_f^\bullet(M)\simeq H_g^\bullet(M)$.
\end{cor}

%%%%%%%%%%%%%%%%%%%%%%%%%%%%%%%%%%
\subsection{Relative cohomology}
%%%%%%%%%%%%%%%%%%%%%%%%%%%%%%%%%%
%We are going to imitate the construction of the de Rham relative cohomology
%done in [BT].\\
%We denote by $\iota$ the inclusion $N\hookrightarrow M$.\\
%For each integer $k$, we put 
%$$\Omega^k(M,N)=\Omega^k(M)\oplus\Omega^{k-1}(N)$$ 
%and we denote by $\tilde{d_f}$ the linear application
%\begin{eqnarray*}
%\tilde{d_f} :
%\Omega^k(M,N) &\longrightarrow& \Omega^{k+1}(M,N)\\
%(\omega,\theta) &\longmapsto& (d_f\omega,\iota^\ast\omega-d_f\theta)
%\end{eqnarray*}
%where $d_f\theta=d_{f_{|N}}\theta$.\\
%It is fairly easy to see that $\tilde{d_f}\circ\tilde{d_f}=0$.\\
%We denote by $H_f^\bullet(M,N)$ the cohomology of this complex. We call it
%{\it relative cohomology} of $N\stackrel{\iota}\hookrightarrow M$ for the
%operator $d_f$.

             %%%%%%%%%%%%%%%%%%%%%%%%%%%  

Let $N$ be a submanifold (eventually with boundary) of $M$. We assume
that $N$ is not included in $S$ and we denote by $\iota$ the inclusion
$N\hookrightarrow M$. The \emph{relative cohomology} groups
$H_f^\bullet(M,N)$ are defined exactly as in the case of the de Rham
theory (see, e.g., the construction done in \cite{BT}).
            
             %%%%%%%%%%%%%%%%%%%%%%%
%If we define $J:\Omega^{k-1}(N)\longrightarrow \Omega^k(M,N)$ and
%$K:\Omega^k(M,N)\longrightarrow \Omega^k(M)$ by 
%$$J(\theta)=(0,\theta)\quad {\mbox { and }}\quad K(\omega,\theta)=\omega\,,$$
%we then have the short exact sequence
%$$0\longrightarrow\Omega^{k-1}(N)\stackrel{J}\longrightarrow\Omega^k(M,N)\stackrel{K}
%\longrightarrow \Omega^k(M)\longrightarrow 0$$
%In the same way as for the de Rham cohomology (see [BT]), since $J$ and $K$ 
%commute with
%$d_f$ et $\tilde{d_f}$ we have the following property.
             %%%%%%%%%%%%%%%%%%%%%

As in case of the de Rham cohomology, we have a long exact sequence
for the pair $(M,N)$:

\begin{thm}
\label{corel}
There is a long exact sequence
\[
\xymatrix{ \cdots\ar[r]& H^{k-1}_f(N)\ar[r]& H^k_f(M,N)\ar[r]&
H^k_f(M)\ar[r]^{\iota^\ast}& H^k_f(N)\ar[r]& \cdots
}
\]
\end{thm}

\begin{cor}
If $M\setminus S$ is dense in $M$, we have $H^0_f(M,N)=\{0\}$ (and
also $H^0_{f,p}(M,N)=\{0\}$).
\end{cor}

\begin{proof}
Apply Proposition \ref{prop:H0cst} and Theorem \ref{corel}.
\end{proof}

Now, assume that $N$ is the closure of an open subset of $M$ instead
of a manifold. We can still define the relative cohomology
$H^\bullet_f(M,N)$. In fact, if we denote by $\Omega^k_N(M)$ the
vector space formed by the $k$-forms which vanish on $N$, then
exterior differentiation $d:\Omega^k_N(M)\to \Omega^{k+1}_N(M)$ is
well defined. Indeed, if $\alpha\in\Omega^k_N(M)$, then $\alpha=0$ on
the interior of $N$, and thus $d\alpha=0$ on $N$, i.e.,
$d\alpha\in\Omega^{k+1}_N(M)$. It follows that the differential
operator $d_f:\Omega^k_N(M)\to \Omega^{k+1}_N(M)$ is also well
defined. Again, imitating the de Rham case, one obtains:

\begin{prop}
If $N$ is the closure of an open subset of $M$ then the 
cohomology of the complex $(\Omega^\bullet_N(M),d_f)$ is isomorphic to
the cohomology $H^\bullet_f(M,N)$.
\end{prop}

\subsection{Excision} We leave it to the reader to check that the
following version of the excision property also holds (again, the
proof is similar to the de Rham case):

\begin{thm}
Let $U$ be an open subset of $M$ with closure in the interior of
$N$. Then, the inclusion $j:(M\setminus U, N\setminus
U)\hookrightarrow (M,N)$ induces an isomorphism
\[ j^*:H^\bullet_f(M,N)\longrightarrow H^\bullet_f(M\setminus
U,N\setminus U).\]
\end{thm}

%%%%%%%%%%%%%%%%%%%%%%%%%%%%%%%%%%%%
\subsection{The Mayer-Vietoris sequence}

Since the differential $d_f$ commutes with the restrictions to open
subsets, one can construct, in the same way as for the de Rham
cohomology (see \cite{BT}), a Mayer-Vietoris exact sequence.

\begin{thm}
\label{MV}
If ${\mathcal U}=(U,V)$ is an open cover of $M$, we have the long
exact sequence 
\[
\xymatrix{
\hdots\ar[r]& H^{k-1}_f(U\cap V)\ar[r] & H^k_f(M) \ar[r]^R & 
H^k_f(U)\oplus H^k_f(V)\ar[r]^J & H^k_f(U\cap V)\ar[r] &\hdots}
\]
where for $[\omega]\in H^k_f(M)$ and $([\alpha],[\beta])\in
H^k_f(U)\oplus H^k_f(V)$, we define 
\[ 
R([\omega])=([\omega|_U],[\omega|_V]),\qquad 
J([\alpha],[\beta])=[\alpha|_{U\cap V}-\beta|_{U\cap V}].
\]
\end{thm} 

%%%%%%%%%%%%%%%%%%%%%%%%%%%%%%%%%%%%%%%
\subsection{Homotopy invariance}
We need to define an appropriate notion of homotopy.
Assuming a more functorial approach, let us think of a pair $(M,f)$ as
an \emph{object}. In order to think of $H^\bullet_f(M)$ as a functor, we need
a notion of \emph{morphism} between such pairs:

\begin{defn}
Let $M$ and $N$ two differentiable manifolds with smooth functions $f$
and $g$, respectively. A \emph{morphism} $\Phi$ from the
pair $(M,f)$ to the pair $(N,g)$ is a pair $(\phi,a)$ formed by a 
smooth map $\phi:M \to N$ and smooth function $a:M\to\R$, such
$a$ does not vanish on $M$ and $g\circ\phi=af$.
\end{defn}

We will say that the pairs $(M,f)$ and $(N,g)$ are \emph{equivalent} if there
exists a morphism $\Phi=(\phi,a)$ between these two pairs where $\phi$ 
is a diffeomorphism. This notion of equivalence between the pairs is
sometimes called ``contact equivalence'' in singularity theory.
 
A morphism $\Phi=(\phi,a)$ from the pair $(M,f)$ to the pair $(N,g)$
induces a chain map 
$\Phi^\ast:(\Omega^\bullet(N),d_g)\to(\Omega^\bullet(M),d_f)$ defined by:
\[
  \Phi^\ast:\Omega^k(N)\to\Omega^k(M),\quad
  \omega\longmapsto\frac{\phi^\ast\omega}{a^k}.
\]
and this map induces an homomorphism in cohomology
$\Phi^\ast: H^\bullet_g(N)\to H^\bullet_f(M)$. If $\Phi$ is an
equivalence this map is an isomorphism.

Now, we come back to our problem:

\begin{defn}
A \emph{homotopy} from the pair $(M,f)$ to the pair $(N,g)$ is given by
two smooth maps
\[ h:M\times [0,1]\to N,\quad a:M\times [0,1]\to \R,\]
such that for each $t\in [0,1]$, we have a morphism 
\[ H_t\equiv(h(\cdot,t),a(\cdot,t)):(M,f)\to (N,g)\] 
(i.e., $a$ does not vanish and $g\circ h(x,t)=a(x,t)f(x)$).
\end{defn}

If $H=(h,a)$ is a homotopy from $(M,f)$ to $(N,g)$, we obtain
%We have linear applications
%\begin{eqnarray*}
%{\rm H}_t^\ast:\Omega^k(N) &\longrightarrow& \Omega^k(M)\\
%              \omega &\longmapsto& \frac{h_t^\ast\omega}{a_t^k}
%\end{eqnarray*}
%These applications commute with $d_f$ and $d_g$. 
a map at the cohomology level 
\[ H_t^\ast:H^\bullet_g(N)\to H^\bullet_f(M).\]

The problem of homotopy invariance is the following: given a homotopy
$H$, from $(M,f)$ to $(N,g)$, is it true that $H_0^\ast=H_1^\ast$ at
the cohomology level?  For general pairs $(M,f)$ and $(N,g)$ this
seems to be a hard problem.  If the complements of the zero level sets
of $f$ and $g$ are dense sets, then in degree zero we do have
$H_0^\ast=H_1^\ast:H^0_f(M)\to H^0_g(N)$. But for higher degree, this is a
much more difficult problem. In the next section, we give some partial
results in the regular case.

\begin{rem}
One can express the notion of homotopy in terms of singular forms. In
fact, it is easy to check that under the correspondence between
singular k-forms $\omega\in\Omega_f^k(M)$ and k-forms
$f^k\omega\in\Omega^k(M)$ (see the the proof of Proposition
\ref{autreexpress}), the map $H^*_t:\Omega^k(N)\to\Omega^k(M)$
corresponds to the pullback $h_t^*:\Omega^k_g(N)\to\Omega^k_f(M)$.
\end{rem}

%%%%%%%%%%%%%%%%%%%%%%%%%%%%%%%%%%%%%%%%%%%%%%%%%%%%%%%%%%%%%%%%%%%%%%%%%%%%%%
\section{The regular case}
\label{section:regular}
In this section, we consider the case where $0$ is a regular value of
$f$. The subset $S=f^{-1}(\{0\})$ in then an embedded submanifold of 
$M$. In order to simplify the exposition we assume that $S$ is connected.

\subsection{Computation of the cohomology.}
It follows from Proposition \ref{prop:H0cst} that $H^0_f(M)=\R$. 
Our main result is the following:

\begin{thm}
\label{thm:regular}
For each $k\geq 1$, there is an isomorphism
\[ H^k_f(M)\simeq H^k_{dR}(M)\oplus H^{k-1}_{dR}(S).\] 
\end{thm}

Before we start the proof we need to introduce some notation.

Let $U\subset U^\prime$ be tubular neighborhoods of $S$. We may assume that 
$U=S\times ]-\epsilon,\epsilon[$ and 
$U^\prime=S\times ]-\epsilon^\prime,\epsilon^\prime[$, with 
$\epsilon^\prime>\epsilon$, and that 
\[
f|_{U^\prime}:S\times ]-\epsilon^\prime,\epsilon^\prime[\to\R, \quad
(x,t)\longmapsto t.
\]
We denote by $\pi$ the projection $U^\prime\to S$.

Let $\rho:\R\to\R$ be a smooth function which is 1 on
$[-\epsilon,\epsilon]$ and has support contained in
$[-\epsilon^\prime,\epsilon^\prime]$. Note that the function
$\rho\circ f$ is 1 on $U$, and we claim that we can assume that the
function $\rho\circ f$ vanishes on $M\setminus U^\prime$. Indeed, let 
$W=\{x\in M\,:\, |f(x)|<\varepsilon^\prime\}$. If
$W=U^\prime$ there is nothing to prove. If not, we have 
$W=U^\prime \cup V$ where $U^\prime$ and $V$ are disjoint open sets.
Then, there exists a smooth function ${\tilde f}$ which equals $f$ on
$U^\prime$ and such that $|f|>\varepsilon^\prime$ on $V$. By
Corollary \ref{cor:deformation}, we can replace $f$ by ${\tilde f}$.

If $\nu$ is a form on $S$, we will denote by $\ove{\nu}$ the form 
$\rho(f)\pi^\ast\nu$. Notice that 
\[
d\ove{\nu}=\rho(f)\pi^\ast(d\nu)+\rho^\prime(f)df\wedge\pi^\ast\nu,
\] 
so we conclude that
\begin{equation}
df\wedge d\ove{\nu}=df\wedge \ove{d\nu}.
\label{eqn:barre}
\end{equation}

\begin{proof}[Proof of Theorem \ref{thm:regular}] 
We split the proof into several lemmas.

\begin{lem}
Any $k$-form $\omega$ on $M$ can be decomposed, in an unique way,
as
\begin{equation}
\omega=f^k\omega_k+f^{k-1}\omega_{k-1}+\hdots +f\omega_1+\omega_0,
\label{eqn:decomposition1}
\end{equation}
where $\omega_i=\ove{\mu_i}+df\wedge\ove{\nu_i}$, with $\mu_i\in
\Omega^k(S)$ and $\nu_i\in \Omega^{k-1}(S)$, for $0\le i\leq k-1$.
\end{lem}

\begin{proof}
Let $\omega$ be a $k$-form on $M$, and write $\omega=(1-\rho(f))\omega
+\rho(f)\omega$. We can decompose $\omega|_{U^\prime}$, in an unique way,
as
\[ 
\omega|_{U^\prime}=f^k\theta_k+f^{k-1}\theta_{k-1}+\hdots+f\theta_1+\theta_0, 
\]
where $\theta_i=\pi^\ast\mu_i+df\wedge\pi^\ast\nu_i$, for $0\le i\leq
k-1$, with $\mu_i\in\Omega^k(S)$ and $\nu_i\in\Omega^{k-1}(S)$. Now,
since $\rho\circ f=0$ on $M\setminus U^\prime$, the $k$-form
$\rho(f)\omega$ may be written as
\[
\rho(f)\omega=f^k\rho(f)\theta_k+f^{k-1}(\ove{\mu_i}+df\wedge\ove{\nu_i})
+\hdots+(\ove{\mu_0}+df\wedge\ove{\nu_0}).
\]
On the other hand, since $1-\rho(f)$ is 0 on a neighborhood of $S$, we can 
write $(1-\rho(f))\omega=f^k\zeta$ for some $k$-form $\zeta$, and the
result follows.
\end{proof}

In the sequel we denote by $\Phi$ the linear application
\[
\Omega^k(M)\oplus\Omega^{k-1}(S)\to \Omega^k(M),\quad
(\alpha,\beta)\longmapsto f^k\alpha+f^{k-1}df\wedge\ove{\beta}.
\]
If $(\alpha,\beta)\in\Omega^k(M)\oplus\Omega^{k-1}(S)$, with $d\alpha=0$
and $d\beta=0$ then, using (\ref{eqn:barre}), we find
\begin{align*}
d_f\big(\Phi(\alpha,\beta)\big)&= f^{k+1}d\alpha-f^kdf\wedge d\ove{\beta},\\
                               &= f^{k+1}d\alpha-f^kdf\wedge \ove{d\beta}=0.
\end{align*}
Similarly, one checks that if $\mu\in\Omega^{k-1}(M)$
and $\nu\in\Omega^{k-2}(S)$, then
\[ \Phi(d\mu,d\nu)=d_f(f^{k-1}\mu-f^{k-2}df\wedge\ove{\nu}).\]
We conclude that $\Phi$ induces a map at the level of cohomology
\begin{align*}
\Phi : H^k_{dR}(M)\oplus H_{dR}^{k-1}(S) &\to H_f^k(M),\\
([\alpha],[\beta]) &\longmapsto [f^k\alpha+f^{k-1}df\wedge\ove{\beta}].
\end{align*}

\begin{lem}
If $k>1$, $\Phi$ is surjective.
\end{lem}

\begin{proof} 
Let $\omega$ be a $k$-form on $M$ with $d_f\omega=0$.  If we decompose
$\omega$ as in (\ref{eqn:decomposition1}), we obtain
\begin{multline*}
d_f\omega=f^{k+1}d\omega_k+f^kd\omega_{k-1}
+f^{k-1}(d\omega_{k-2}-df\wedge\omega_{k-1})\\
+\cdots+f(d\omega_0-(k-1)df\wedge\omega_1)-kdf\wedge\omega_0=0.
\end{multline*}
If we restrict to $U$, we get by uniqueness of the decomposition 
$df\wedge{\omega_0}|_{U}=0$, i.e. $df\wedge\pi^\ast\mu_0=0$ and so, $\mu_0=0$.
We conclude that $\omega_0=df\wedge\ove{\nu_0}$.

Now set $\gamma_0\equiv\frac{1}{k-1}\ove{\nu_0}$. We have 
\[
\omega+d_f\gamma_0=f^k\omega_f+f^{k-1}\omega_{k-1}+\hdots+f^2\omega_2
+f(\omega_1+d\gamma_0).
\]
Noting that $d\gamma_0=\frac{1}{k-1}\ove{d\nu_0}
+\frac{\rho^\prime(f)}{k-1}df\wedge\pi^\ast\nu_0$, writing
$d_f(\omega+d_f\gamma_0)=0$ and restricting to $U$, we obtain 
$\mu_1+\frac{1}{k-1}d\nu_0=0$. Therefore:
\[ 
\omega_1+d\gamma_0=df\wedge(\ove{\nu_1}
+\frac{\rho^\prime(f)}{k-1}\pi^\ast\nu_0).
\]
Thus, if we put 
$\gamma_1=\frac{1}{k-2}\ove{\nu_1}+\frac{\rho^\prime(f)}{k-1}\pi^\ast\nu_0$,
we get 
\[
\omega+d_f(f\gamma_1+\gamma_0)=f^k\omega_f+f^{k-1}\omega_{k-1}+\hdots
+f^3\omega_3+f^2(\omega_2+d\gamma_1).
\]
This way, we can construct $\gamma_0,\gamma_1,\hdots,\gamma_{k-2}$, with
\[
\gamma_{k-2}=\ove{\nu_{k-2}}+\frac{\rho^\prime(f)}{2}\pi^\ast\nu_{k-3}+
\hdots+\frac{\rho^{(k-2)}(f)}{(k-1)!}\pi^\ast\nu_0,
\]
such that
\[
\omega+d_f(\gamma_0+\hdots+f^{k-2}\gamma_{k-2})=f^k\omega_k+
f^{k-1}(\omega_{k-1}+d\gamma_{k-2}).
\]

Now, writing $d_f(\omega+d_f(\gamma_0+\hdots+f^{k-2}\gamma_{k-2}))=0$ and
restricting to $U$, we obtain $\mu_{k-1}=-d\nu_{k-2}$ (using $\rho\circ f=1$)
and $d{\omega_{k-1}}_{|U}=0$. Consequently, we have 
\[
\omega_{k-1}+d\gamma_{k-2}=df\wedge\ove{\nu_{k-1}}+\eta,
\]
where 
\[
\eta=df\wedge [\rho^\prime(f)\pi^\ast\nu_{k-2}+\hdots
+\frac{\rho^{(k-1)}(f)}{(k-1)!}\pi^\ast\nu_0]+ 
\frac{\rho^\prime(f)}{2}\pi^\ast d\nu_{k-3}+
\hdots+\frac{\rho^{(k-2)}(f)}{(k-1)!}\pi^\ast d\nu_0.
\]
Since $d{\omega_{k-1}}|_{U}=0$ and ${\eta}|_{U}=0$, we obtain
$d\nu_{k-2}=0$. On the other hand, since $\eta$ is zero on a
neighborhood of $S$, we can write $\eta=f\xi$. We conclude that
\[
\omega=f^k(\omega_k+\xi)+f^{k-1}df\wedge\ove{\nu_{k-2}}+d_f\gamma,
\]
where $\gamma=\gamma_0+\hdots+f^{k-2}\gamma_{k-2}$. We have seen, that 
$d\nu_{k-2}=0$. Now, writing $d_f\omega=0$, we see that 
$d(\omega_k+\xi)=0$. This shows that $\omega$ is in the image of $\Phi$.
\end{proof}

\begin{lem}
If $k>1$, $\Phi$ is injective.
\end{lem}

\begin{proof} 
Let $(\alpha,\beta)$ in 
$\Omega^k(M)\oplus\Omega^{k-1}(S)$ with $d\alpha=0$ and $d\beta=0$. We assume
that $f^k\alpha+f^{k-1}df\wedge\ove{\beta}=d_f\gamma$, where 
$\gamma\in\Omega^{k-1}(M)$.

We decompose $\gamma$ as in (\ref{eqn:decomposition1}), i.e.,
\[
\gamma=f^{k-1}\gamma_{k-1}+f^{k-2}\gamma_{k-2}+\hdots+f\gamma_1+\gamma_0,
\]
with, for $i\leq k-2$, $\gamma_i=\ove{\mu_i}+df\wedge\ove{\nu_i}$, $\mu_i$ and
$\nu_i$ are forms on $S$. We have 
\begin{multline*}
d_f\gamma=f^kd\gamma_{k-1}+f^{k-1}d\gamma_{k-2}
+f^{k-2}(d\gamma_{k-3}-df\wedge\gamma_{k-2})\\
+\hdots+f(d\gamma_0-(k-2)df\wedge\gamma_1)-(k-1)df\wedge\gamma_0.
\end{multline*}
Restricting to $U$, we obtain
\begin{align}
                  df\wedge{\gamma_0}|_{U} &=0       \notag\\
   (d\gamma_0-(k-2)df\wedge\gamma_1)|_{U} &=0       \notag\\
                                          &\vdots   \label{eq:recursion}\\
(d\gamma_{k-3}-df\wedge\gamma_{k-2})|_{U} &=0       \notag\\
                    {d\gamma_{k-2}}|_{U}  &= df\wedge\pi^\ast\beta\notag
\end{align}
The first relation gives $df\wedge\pi^\ast\mu_0=0$ and so, $\mu_0=0$. 
This implies that $\gamma_0=df\wedge{\ove \nu_0}$. Using the second
relation, we then get  
\[ df\wedge \pi^\ast d\nu_0+(k-2)df\wedge\pi^\ast\mu_1=0,\]
which implies $\mu_1=-\frac{1}{k-2}d\nu_0$. In this way, we obtain for
each $i\leq k-2$, 
\[
\mu_i=-\frac{1}{k-1-i}d\nu_{i-1}.
\]

Now, since $\gamma_{k-2}=\ove{\mu_{k-2}}+df\wedge\ove{\nu_{k-2}}$, the
one before the last relation in (\ref{eq:recursion}) gives 
$-df\wedge\pi^\ast d\nu_{k-2}=df\wedge\pi^\ast\beta$, which implies 
$\beta=-d\nu_{k-2}$, i.e., $\beta$ is exact.

On the other hand, we have, for each $1\leq i\leq k-2$,
\begin{align*}
d\gamma_{i-1}-(k-1-i)df\wedge\gamma_i &= d\ove{\mu_{i-1}}-df\wedge 
d\ove{\nu_{i-1}}-(k-1-i)df\wedge\ove{\mu_i},\\
                                      &= \ove{d\mu_{i-1}}
+\rho^\prime(f)df\wedge\pi^\ast \mu_{i-1}\\
&\qquad \qquad \qquad -df\wedge[\ove{d\nu_{i-1}}+
(k-1-i)\ove{\mu_i}],\\
                               &= -\frac{\rho^\prime(f)}{k-1-i} df\wedge
\pi^\ast d\nu_{i-1},
\end{align*}
and
\[
d\gamma_{k-2}=df\wedge\ove{\beta}+\rho^\prime(f) df\wedge\pi^\ast\mu_{k-2}=
df\wedge\ove{\beta}-\rho^\prime(f) df\wedge\pi^\ast d\nu_{k-3}.
\]
We conclude that 
\[
f^k\alpha=f^k d\big(\gamma_{k-1}+
\frac{\rho^\prime(f)}{f}df\wedge\pi^\ast\nu_{k-3}+
\frac{\rho^\prime(f)}{f^2}df\wedge\pi^\ast\nu_{k-3}+\hdots+
\frac{\rho^\prime(f)}{(k-2)f^{k-1}}df\wedge\pi^\ast\nu_{0}\big).\]
Therefore, $\alpha$ is exact.
\end{proof}

This shows that $\Phi$ is bijective for $k>1$. On the other hand, we
have

\begin{lem}
If $k=1$, $\Phi$ is bijective.
\end{lem}

\begin{proof}
To prove that $\Phi$ is surjective, let $\omega$ be a 1-form on $M$
with $d_f\omega=0$. We write $\omega=f\omega_1+\omega_0$ with 
$\omega_0=\ove{\mu_0}+df\wedge\ove{\nu_0}$ ($\mu_0\in\Omega^1(S)$,
$\nu_0\in\Omega^0(S)$).
We write $d_f\omega=0$ an we restrict to $U$. We obtain $\mu_0=0$ hence,
$\omega_0=df\wedge\ove{\nu_0}$. Moreover, we have $d{\omega_0}_{|U}=0$ which gives 
$d\nu_0=0$. It follows that $d\omega_1=0$.

Now to prove that $\Phi$ is injective, let $\alpha\in \Omega^1(S)$ and
$\beta\in\Omega^0(S)$ with $d\alpha=0$ and $d\beta=0$. We suppose that 
$$f\alpha+df\wedge\ove{\beta}=d_f\gamma=fd\gamma$$
where $\gamma\in\Omega^0(M)$. Restricting to $S$, we obtain $\beta=0$. This implies
$\alpha=d\gamma$.
\end{proof}

We have establish that $\Phi$ is an isomorphism for all $k\ge 1$ so
Theorem \ref{thm:regular} follows.
\end{proof}

\begin{rem}
\label{rk:cohomelrose}
Comparing this result with Proposition 2.49 in \cite{Me}, we see that
the cohomology of the Lie algebroid attached to a function constructed
in Section \ref{sec:Lie algebroid} is isomorphic to the cohomology of
the Melrose Lie algebroid.
\end{rem}

\begin{rem}
\label{rk:warning0}
For $k-p>0$ it is possible to adapt this proof in order to compute the
cohomology $H_{f,p}^\bullet(M)$. For $k-p<0$ the decomposition
$(\ref{eqn:decomposition1})$ is no longer valid. For $k=p$ the
expression for $H^p_{f,p}(M)$ is not so nice. For instance, if $p=1$,
and if $H^1_{dR}(M)=\{0\}$, we can show that the space $H^1_{f,1}(M)$
has infinite dimension. In fact, the space $Z^1_{f,1}(M)$ of 1-cocycles is 
$\{dh\,|\,h\in C^\infty(M)\}$ which is isomorphic, via
exterior differentiation $d$, to the space $C^\infty_0(M)$ of
functions which vanish in at least a point of $M$. Similarly,
the space $B^1_{f,1}(M)$ of 1-cobords is isomorphic, via $d$, to the
ideal of $C^\infty_0(M)$ spanned by $f$. Therefore, the quotient
$C^\infty_0(M)/(f)$ has infinite dimension.
\end{rem}

\begin{example}
Let $M=\{(x_1,\hdots,x_n)\in \R^n:x_1^2+\cdots+x_n^2<2\}$ be an open
ball and $f:M\to\R$ the function
$f(x_1,\hdots,x_n)=x_1^2+\cdots+x_n^2-1$, so that $S\subset M$ is the
$(n-1)$-sphere. Then, 
\begin{align*}
H^0_f(M)&= H^1_f(M)= H^n_f(M)= \R,\\
H^k_f(M)&=\{0\},\quad \text{ if } 2\leq k\leq n-1.
\end{align*}
\end{example}

\begin{example}
Let $M=\{(x_1,\hdots,x_{n+1})\in \R^n:x_1^2+\cdots+x_{n+1}^2=1\}$ be
the $n$-sphere and $f:M\to\R$ the function $f(x_1,\hdots,x_{n+1})=x_1$,
so that $S$ is the equator. Then,
\begin{align*}
H^k_f(M) &= \R\quad \text{ if } k=0,1,\\
H^k_f(M) &= \{0\}  \quad \text{ if } 2\leq k\leq n-1,\\ 
H^n_f(M) &= \R^2.
\end{align*}
\end{example}

\begin{example} (Poisson geometry) 
Recall the identification explained in Section \ref{sec:2 dim} between
the cohomology $H_f^\bullet(M)$ and Poisson cohomology in dimension 2.
It leads immediately to the following result, which generalizes a
result due to Radko \cite{R} for the compact case:

\begin{thm}
\label{thm:Poisson cohomology 2}
Let $(M,\Pi)$ be an orientable 2-dimensional Poisson manifold with
singular set $S$. Assume that the contraction of the Poisson tensor 
$\Pi$ with a volume form on $M$ is a regular function in a neighborhood
of $S$. Then the Poisson cohomology of $(M,\Pi)$ is
\[ H^k_\Pi(M)\simeq H^k_{dR}(M)\oplus H^{k-1}_{dR}(S).\]
\end{thm}
\end{example}

\subsection{Homotopy invariance in the regular case.}
In the regular case we are able to prove homotopy invariance:

\begin{prop}
\label{prop:homotopy1}
Let $U$ and $W$ be tubular neighborhoods of $S_f=f^{-1}(0)$
and $S_g= g^{-1}(0)$, respectively. We assume that $f$ and $g$ do
not have singularities on $U$ and $W$. If $H_t$ is a homotopy from
$(U,f)$ to $(W,g)$. Then the induced linear applications between the
cohomology spaces are the same: $H^\ast_1=H^\ast_0$.
\end{prop}

\begin{proof}
We can assume that $U=S_f\times ]-\varepsilon,\varepsilon[$ and 
$W=S_g\times ]-\varepsilon^\prime,\varepsilon^\prime[$, with
\[ 
(x,\rho)\stackrel{f}{\longmapsto} \rho\quad \text{ and }
\quad (y,\tau)\stackrel{g}{\longmapsto} \tau.
\]
By Proposition \ref{prop:H0cst} we can take $k\geq 1$. We denote by
$\Psi_f$ and $\Psi_g$ the linear maps:
\begin{align*}
\Psi_f: H^k_{dR}(U)\oplus H^{k-1}_{dR}(U)&\to H^k_f(U)\\
([\alpha],[\beta])&\longmapsto [\rho^k\alpha+\rho^{k-1}d\rho\wedge\beta],\\
\Psi_g: H^k_{dR}(W)\oplus H^{k-1}_{dR}(W)&\to H^k_f(W)\\
([\alpha],[\beta])&\longmapsto [\tau^k\alpha+\tau^{k-1}d\tau\wedge\beta],
\end{align*}
which, by Theorem \ref{thm:regular}, are isomorphisms. 

Now, we set $K_t^\ast=\Psi_f^{-1}\circ H_t^\ast\circ\Psi_g$, for every
$t\in [0,1]$. If $([\alpha],[\beta])\in H^k_{dR}(W)\oplus H^{k-1}_{dR}(W)$, 
we have
\begin{align*}
H_t^\ast\big(\Psi_g([\alpha],[\beta])\big)&= 
\big[\frac{h_t^\ast(\tau^k\alpha+\tau^{k-1}d\tau\wedge\beta)}{a_t^k}\big],\\
 &=\big[\frac{a_t^k\rho^kh_t^\ast\alpha+a_t^{k-1}\rho^{k-1}(\rho
   da_t\wedge h_t^\ast\beta+a_td\rho\wedge h_t^\ast\beta)}{a_t^k}\big],\\
 &=[\rho^kh_t^\ast\alpha+\rho^k\frac{da_t}{a_t}\wedge h_t^\ast\beta
   +\rho^{k-1}d\rho\wedge h_t^\ast\beta],\\
 &=[\rho^kh_t^\ast\alpha+\rho^{k-1}d\rho\wedge h_t^\ast\beta
   +\rho^kd(log|a_t|h_t^\ast\beta)].
\end{align*}
We conclude that
\begin{align*}
K_t^\ast\big([\alpha],[\beta]\big)
   &=\big([h_t^\ast\alpha+d(log|a_t|h_t^\ast\beta)],[h_t^\ast\beta]\big),\\
   &=\big([h_t^\ast\alpha],[h_t^\ast\beta]\big).
\end{align*}
Since the de Rham cohomology is homotopy invariant, we have
$K_1^\ast=K_0^\ast$ and it follows that $H_1^\ast=H_0^\ast$.
\end{proof}

\begin{prop}
\label{prop:homotopy2}
Let $H_t$ be a homotopy from $(M,f)$ to $(N,g)$. We assume that $f$ and 
$g$ are regular on tubular neighborhoods of $S_f$ and $S_g$. 
If $H^{k-1}_{dR}(S)$ is trivial, then the linear maps $H_0^\ast$ and
$H_1^\ast$ from $H^k_g(N)$ to $H^k_f(M)$ coincide.
\end{prop}

\begin{proof}
Note that the assumptions imply that $k\geq 2$.

Let $U$ and $W$ be tubular neighborhoods of $S_f$ and $S_g$ such that
$f$ and $g$ are regular on these neighborhoods. We can assume that H sends
$W$ onto $U$, and we set $V=M\setminus S_f$ and $Z=N\setminus S_g$.

Let $\omega$ be in $Z^k_g(N)$. According to the previous proposition, we have
\[
(H_1^\ast\omega)|_{U}=
(H_0^\ast\omega)|_{U}+d_f\alpha_U,\quad \alpha_U\in\Omega^{k-1}(U).
\]
On the other hand, since $f$ and $g$ do not vanish on $V$ and $Z$ 
and since the de Rham cohomology is homotopy invariant, we have
\[
(H_1^\ast\omega)|_{V}=
(H_0^\ast\omega)|_{V}+d_f\alpha_V,\quad \alpha_V\in\Omega^{k-1}(V).
\] 
Therefore, we obtain
\[
d_f({\alpha_U}|_{U\cap V}-{\alpha_V}|_{U\cap V})=0,
\]
i.e., ${\alpha_U}|_{U\cap V}-{\alpha_V}|_{U\cap V}\in Z^{k-1}_f(U\cap
V)$.

Now, since $H^{k-1}_f(U\cap V)\simeq H^{k-1}_{dR}(U\cap V)
\simeq \big(H^{k-1}_{dR}(S)\big)^2=\{0\}$, there exists 
$\beta_{U\cap V}\in\Omega^{k-2}(U\cap V)$ such that 
\[
{\alpha_U}|_{U\cap V}-{\alpha_V}|_{U\cap V}=d_f\beta_{U\cap V}.
\]
{From} the exactness of the Mayer-Vietoris short exact sequence for
de Rham cohomology, there exist $\alpha^\prime_U\in\Omega^{k-2}(U)$ 
and $\alpha^\prime_V\in\Omega^{k-2}(V)$ such that
$\beta_{U\cap V}={\alpha^\prime_V}|_{U\cap V}-{\alpha^\prime_U}|_{U\cap V}$.
It follows that
\[
(\alpha_U+d_f\alpha^\prime_U)|_{U\cap V}=
(\alpha_V+d_f\alpha^\prime_V)|_{U\cap V}.
\]
Hence, there exists $\eta\in\Omega^{k-1}(M)$ such that 
\[
\eta|_{U}=\alpha_U+d_f\alpha^\prime_U \quad \text{ and }\quad
\eta|_{V}=\alpha_V+d_f\alpha^\prime_V.
\]
This gives
\[
(d_f \eta)|_{U}=(H_1^\ast\omega-H_0^\ast\omega)|_{U}
\quad \text{ and }\quad
(d_f \eta)|_{V}=(H_1^\ast\omega-H_0^\ast\omega)|_{V},
\]
which shows that
\[ 
H_1^\ast\omega-H_0^\ast\omega=d_f\eta.
\]
\end{proof}

%%%%%%%%%%%%%%%%%%%%%%%%%%%%%%%%%%%%%%%%%%%%%%%%%%%%%%%%%%%%%%%%%%%%%%%%%%%%%%%

\section{One step to the complex case}
The definition of the cohomology $H^\bullet_f(M)$ readily extends to
complex manifolds. In this section we study the local case and give an
application of this cohomology to the study of the topology of the
complement of an hypersurface.

We feel that this cohomology may have others applications in algebraic
geometry or in analytic geometry, and that from it one may be able to
obtain more information on the topology of the complement of the zeros
of a function $f$.

\subsection{Cohomology in the the local case.}
In this paragraph we give an overview of the results we have found in
\cite{M2,M3}. There we consider a germified version of the cohomology:
we let $\Omega^k(\C^n)$ denote the space of germs at 0 of analytic
$k$-forms, and we let $H_{f,p}^\bullet(\C^n)$ denote the cohomology of
the chain complex $(\Omega^k(\C^n),d_f^{(p)})$. We consider only the
groups $H_{f,p}^{n-1}(\C^n)$ and $H_{f,p}^{n}(\C^n)$. The other groups
are usually trivial, with the exception of $H^0$ and $H^1$ (see
\cite{M2,M3}). 

We will assume that the function $f$ is a quasi-homogeneous
polynomial on $\C^n$ of degree $N$, with respect to the weights
$w_1,\dots,w_n$, and with an isolated singularity at 0. We denote by
$c$ the \emph{Milnor number} of the singularity, i.e., the dimension
of the vector space $Q_f={\mathcal O}_n/I_f$ where ${\mathcal O}_n$ is
the space of germs of analytic functions and $I_f$ the ideal spanned
by the first derivatives of $f$. Also, for every positive integer $q$,
we denote by $h^{q,n-q}$ the dimension of ${(Q_f)}_{qN-w_1-\cdots -w_n}$, 
the quasi-homogeneous part of degree $qN-w_1-\hdots -w_n$ of the graded
space $Q_f$. These numbers are the mixed Hodge numbers of the
quasi-homogeneous singularity $f$.

Table 1 summarizes the results obtained in \cite{M2}.
\vskip 15 pt

\begin{center}
\begin{tabular}{|c|c|c|}
\hline
& & \\
      & $\dim H^{n-1}_{f,p}(\C^n)$ & $\dim H^n_{f,p}(\C^n)$\\
& & \\
\hline
& &\\
$0\le p\le n-3$  & $\sum_{i=1}^{n-p-1} h^{i,n-i}$ &
      $c+\sum_{i=1}^{n-p-1} h^{i,n-i}$\\
& &\\
$p=n-2$ & $\infty$                           & $c+h^{1,n-1}$\\
& &\\
$p=n-1$ & ?                                 & $\infty$\\
& &\\
$p\geq n$ & 0 & $c$\\
\hline
\end{tabular}
\vskip 10 pt

\textbf{Table 1}
\end{center}

\vskip 20 pt

\begin{rem}
For $k>0$ denote by $\Omega^k_{rel}(\C^n,f)$ the quotient
$\Omega^k(\C^n)/df\wedge\Omega^{k-1}(\C^n)$. It is easy to check that
the de Rham differential $d$ passes to the quotient, so we get a
complex $(\Omega^\bullet_{rel}(\C^n,f),d)$. The cohomology of this
complex is the well-known \emph{relative cohomology} of the singularity $f$.
This cohomology seems to be linked with the cohomology
$H_{f,p}^\bullet(\C^n)$, but they do not coincide (e.g., compare the
table above with the results in \cite{Ma}). Nevertheless, the
computation of the cohomology $H_{f,p}^\bullet(\C^n)$ presented in
\cite{M2}, uses the vanishing of certain relative cohomology spaces of
$f$.
\end{rem}

\subsection{Cohomology of the complement of an hypersurface}

We shall now explain a method, using the cohomology
$H_{f,p}^\bullet(\C^n)$, to obtain information on the
cohomology of the complement of an hypersurface. More precisely, we
apply this cohomology to determine at which stage a certain spectral sequence 
converging to the cohomology of an hypersurface singularity degenerates.
We then use this to give a new proof of a well-known result of
A.~Dimca (\cite{D1}).

\subsubsection{Local case}
Let $B$ be a small open ball at the origin of $\C^n$. We consider a 
hypersurface singularity $V\subset B$ at the origin. Let $f=0$
be an equation for $V$ in $B$ and denote by $U=B\setminus V$ the
complement. A well-known result of Grothendieck (\cite{Gk}) states
that the cohomology $H^\bullet(U,\C)$ is isomorphic to the cohomology 
of the complex $A^\bullet_0$ of meromorphic differential forms on $B$
with polar singularity along $V$. 

An element $\omega\in A^k_0$ can be written in the form
$\omega=\frac{\alpha}{f^s}$ where $\alpha$ is a holomorphic $k$-form
on $B$. We consider the decreasing filtration:
\[
F^s A^j_0=\left\{
\begin{array}{l}
\left\{\frac{\alpha}{f^{j-s}}:\alpha \text{ holomorphic on }B\right\}
\quad {\text{ if }} j-s\geq 0,\\
\\
\{0\}\quad \text{ if } j-s<0.
\end{array}
\right.
\]
This filtration is exhaustive and bounded above so, it induces a
spectral sequence $\big( E_r(V),d_r\big)$ converging to
$H^\bullet(U,\C)$ (see \cite{Mc}). It is known (see \cite{D2}) that
this spectral sequence degenerates after a finite number of steps. The
problem is to determine this number.

For every $p,q,r$ we set
\[ E_r^{p,q}(V)=\frac{Z_r^{p,q}(V)}{Z_{r-1}^{p+1,q-1}(V)+B_{r-1}^{p,q}(V)}\]
where $Z_r^{p,q}(V)$ and $B_r^{p,q}(V)$ are well-known spaces (see
\cite{Mc}) and $d_r^{p,q}:E_r^{p,q}(V)\to E_r^{p+r,q-r+1}(V)$.
This spectral sequence degenerates at the step $r$, i.e., $E_r=E_\infty$, if 
$d_r^{.,.}=0$. In order to show that $d_r^{.,.}=0$ (for some $r$) it
is sufficient to show  that for every $p,q$ we have 
\begin{equation}
\label{eq:degenerate}
d\big(Z_r^{p,q}(V)\big)\subset B_{r-1}^{p+r,q-r+1}(V).
\end{equation}
If we remark that, for a holomorphic $(p+q)$-form $\alpha$, we have 
$d\big(\frac{\alpha}{f^q}\big)=\frac{d_f^{(p)}\alpha}{f^{q+1}}$,
then we can rewrite (\ref{eq:degenerate}) as:
\begin{itemize}
\item if $\alpha$ is a holomorphic $(p+q)$-form such that $f^r$ 
divides $d_f^{(p)}\alpha$ then, for some holomorphic $(p+q)$-form
$\zeta$, one has $d_f^{(p)}\alpha=d_f^{(p)}(f\zeta)$.
\end{itemize}
\vskip 10 pt

It is known that when the function $f$ is regular one has $E_1=E_\infty$. Now,
we assume that $f$ has an isolated singularity at 0.  In this case, it
is known (see \cite{D2}) that $d_1^{p,q}=0$ if $p+q<n-1$. Let us look
then at $d_r^{p,q}$ with $p+q=n-1$. We assume further that $f$ is a
$W$-quasi-homogeneous polynomial of degree N, where
$W=w_1x_1\frac{\partial }{\partial x_1}+\hdots w_nx_n\frac{\partial
}{\partial x_n}$, with each $w_i$ a positive integer. This means that:
\[ W\cdot f=N f.\]
In \cite{M2}, we have computed the spaces $H^n_{f,p}(B)$ under these
assumptions, and we recall here our results. We set $Q_f={\mathcal
H}(B)/I_f$, where ${\mathcal H}(B)$ is the algebra of holomorphic
functions on $B$ and $I_f$ the ideal spanned by $\frac{\partial
f}{\partial x_1},\dots,\frac{\partial f}{\partial x_n}$.  This vector
space has finite dimension (the Milnor number of $f$) and we let
${\mathcal B}$ denote a monomial basis (for the existence of
such a basis, see \cite{AGV}). Finally, we set
$\nu=dx_1\wedge\hdots\wedge dx_n$.

\begin{thm}[\cite{M2}]
\label{Hn}
Assume that $p<n-1$ and let $\eta\in\Omega^n(B)$. There exist unique
polynomials $h_1,\hdots,h_{n-p}$ (possibly zero) such that:
\begin{enumerate}
\item[(a)] $h_1$ is quasi-homogeneous of degree $N-\sum w_i$;
\item[(b)] $h_j$ for $2\leq j\leq n-p-1$ is a linear combination of monomials 
of ${\mathcal B}$ of degree $jN-\sum w_i$;
\item[(c)] $h_q$ is a linear combination of monomials of ${\mathcal B}$ and
\[
\eta=(h_{n-p}+fh_{n-p-1}+\cdots +f^{n-p}h_1)\nu\pmod{B^n_{f,p}(B)}.
\]
\end{enumerate}
\end{thm}

This theorem allows us to give a new proof of the following result
(see \cite{D1}). 

\begin{cor}
If $f$ is a quasi-homogeneous polynomial with an isolated singularity at 0, 
then the spectral sequence degenerates after the second step, i.e., 
$E_2=E_\infty$.
\end{cor}

\begin{proof}
We only need to consider $d_2^{p,q}$ with $p+q=n-1$. Also, if $q=0$ it
is easy to see that $d^{n-1,0}_2=0$, so we assume $q>0$, i.e., $p<n-1$.

Let $\alpha$ be an $(n-1)$-form on $B$. We will show that if for some $n$-form
$\theta$ on $B$ one has $f^2\theta=d_f^{(p)}\alpha$, then there exists an
$(n-1)$-form $\zeta$ such that $d_f^{(p)}\alpha=d_f^{(p)}(f\zeta)$.

By Theorem \ref{Hn}, if $\theta$ is a holomorphic $n$-form on $B$, 
we have 
\[ 
f\theta=(h_{n-p-1}+fh_{n-p-2}+\cdots +f^{n-p-1}h_1)\nu+d_f^{(p+1)}\zeta,
\]
where $\zeta$ is a holomorphic $(n-1)$-form and the $h_i$ are as in
the theorem. It follows from Lemma 3.12 in \cite{M2} that
$d_f^{(p+1)}\zeta\in I_f$. Since $f$ is also in $I_f$ we must have
$h_{n-p-1}=0$. Since $fd_f^{(p+1)}\zeta=d_f^{(p)}(f\zeta)$, we see that
\[ 
f^2\theta=(f^2h_{n-p-2}+\cdots+f^{n-p}h_1)\nu+d_f^{(p)}(f\zeta).
\]
Hence, if $f^2\theta=d_f^{(p)}\alpha\in B^n_{f,p}(B)$,
we have $(f^2h_{n-p-2}+\hdots+f^{n-p}h_1)\nu\in B^n_{f,p}(B)$. The
previous relation then implies that $h_{n-p-2}=\cdots =h_1=0$.
Therefore, we conclude that $d_f^{(p)}\alpha=f^2\theta=d_f^{(p)}(f\zeta)$.
\end{proof}

%%%%%%%%%%%%%%%%%%%%%%%%%%%%%%%%%%%%%%%%%%%%
\subsubsection{Global projective case}
Let $W$ be the vector field $w_1x_1\frac{\partial }{\partial
x_1}+\hdots w_nx_n\frac{\partial }{\partial x_n}$, with
$w_1,\hdots,w_n$ positive integers. We denote by $\Proj^n(W)$ the
weighted projective space associated to $W$ (see \cite{Dg}). We
consider a quasi-homogeneous polynomial (with respect to $W$)
$f\in\C[x_0,\hdots,x_n]$ of degree $N$ and we denote by $V$ the
hypersurface in $\Proj^n(W)$ with equation $f=0$. Again, the
cohomology $H^\bullet(U,\C)$, where $U=\Proj^n(W)\setminus V$), is
isomorphic to the cohomology of the complex $A^\bullet$ of algebraic
differential forms (see \cite{Gk}). 

An element $\omega$ in $A^k$ can be written as $\omega=\frac{\alpha}{f^s}$ 
where $\alpha$ is a quasi-homogeneous $k$-form of degree $sN$. This
means that $i_W\alpha=0$ and the Lie derivative satisfies ${\mathcal
L}_W\alpha=(sN)\alpha$. In this case we can consider the following
decreasing filtration
\[
{\tilde F}^s A^j=
\left\{
\begin{array}{l}
\left\{
\begin{array}{c}
\frac{\alpha}{f^{j-s}}:\alpha \text{ quasi-homogeneous
of degree } (j-s)N\\ \text{ and } i_W\alpha=0
\end{array}
\right\}\quad \text{ if }
j-s>0,\\
\\
\{0\}\quad \text{ if } j-s\leq 0.
\end{array}
\right.
\]

This filtration induces a spectral sequence $({\tilde E}_r(V),{\tilde
d}_r)$ converging to $H^\bullet(U,\C)$.  If $f$ is regular, this
spectral sequence degenerates after the first step (see
\cite{G}). Now, we assume that $f$ has an isolated singularity at
0. In this case, one knows that ${\tilde d}_1^{p,q}=0$ if $p+q<n-1$
(see \cite{D2}).

\begin{prop}
If $f$ is a quasi-homogeneous polynomial with an isolated singularity at 0, 
then the spectral sequence degenerates after the second step, i.e., 
${\tilde E}_2={\tilde E}_\infty$.
\end{prop}

\begin{proof}
According to \cite{D2}, we only need to consider ${\tilde d}_2^{p,q}$
with $p+q=n-1$. As for the local case, we need to show that if
$\alpha$ is a quasi-homogeneous $(n-1)$-form on $\C^{n+1}$, of degree
$qN$ ($q=n-1-p$) such that $i_W\alpha=0$ and $f^2$ divides
$d_f^{(p)}\alpha$, then there exists a quasi-homogeneous $(n-1)$-form
$\zeta$ which satisfies $i_W\zeta=0$ and $d_f^{(p)}\alpha=d_f^{(p)}(f\zeta)$.

Let us denote by $\eta$ the $n$-form $d_f^{(p)}\alpha$. It is easy to
check that $i_W\eta=0$. Therefore, we have $\eta=i_W(g\nu)$, where $g$ is  
some quasi-homogeneous polynomial of degree $(q+1)N-\sum w_i$. Set
$\sigma=i_W\nu$, so that $\eta=g\sigma$. By Lemma 3.26 in \cite{M2},
we have 
\[
\eta\in B^n_{f,p}(\C^{n+1}) \Longleftrightarrow 
g\nu\in B^{n+1}_{f,p+1}(\C^{n+1}).
\]
Since $f^2$ divides $g$, we can write $g\nu=f^2\xi$, where $\xi$ is some 
quasi-homogeneous $(n+1)$-form on $\C^{n+1}$. We then have $f^2\xi\in
B^{n+1}_{f,p+1}(\C^{n+1})$. Now, it is possible to adapt
the argument we gave above in the local case (theorem \ref{Hn} is still valid
in a polynomial version because of the homogeneity of the operators
$d_f^{(p)}$) to obtain $f^2\xi=d_f^{(p+1)}(f\mu)$, where $\mu$ is a
quasi-homogeneous $n$-form. We conclude that
\[ \eta=i_W\big(d_f^{(p+1)}(f\mu)\big)=-d_f^{(p)}\big(f(i_W\mu)\big).\]
\end{proof}

\end{document}